\documentclass{article-hermes}

\usepackage{graphicx}

\usepackage[bookmarks=true, bookmarksnumbered=true, bookmarksopen=true,
		unicode=true, colorlinks=true,
		linkcolor=blue,citecolor=blue,filecolor=blue,urlcolor=blue,
		pdfstartview=FitH
		]{hyperref}

\begin{document}

\title[Propriétés petit-monde et sans-échelle]%
      {Étude de l'omniprésence des propriétés petit-monde et sans-échelle} 

\author{Vincent Labatut}

\address{
Université Galatasaray, Ç{\i}ra\u{g}an cad. n°36, Ortaköy 34357, \.{I}stanbul, Turquie\\
LIA, Université d'Avignon, Agroparc BP 1228, 84911 Avignon Cedex 9, France\\[3pt]
\href{mailto:vincent.labatut@univ-avignon.fr}{vincent.labatut@univ-avignon.fr}}

\resume{Les propriétés petit-monde et sans-échelle ont été identifiées dans les réseaux complexes réels à la fin des années 90. Leur étude a permis de faire progresser la compréhension de certains systèmes, et elles ont fait l'objet de très nombreux travaux. Peut-être pour cette raison, on rencontre très fréquemment dans la littérature l'idée qu'il s'agirait de propriétés quasi-universelles, valables pour presque tous les réseaux complexes. Or, à notre connaissance, aucune étude à grande échelle n'a été menée pour véritablement répondre à cette question. Dans ce travail, nous tirons parti, d'une part, de la grande quantité de données qui est maintenant disponible en ligne pour constituer une large collection de réseaux complexes, et d'autre part d'outils d'analyse définis récemment, pour éprouver objectivement cette hypothèse d'omniprésence. Il ressort que si une large majorité des réseaux traités est bien petit-monde, ce n'est en revanche pas le cas pour la propriété sans-échelle, puisque le degré d'une proportion significative de ces réseaux n'est pas distribué selon une loi de puissance.}

\abstract{The small-world and scale-free properties were identified in real-world complex networks at the end of the 90s. Their analysis led to a better understanding of the dynamics and functioning of certain systems, and they were studied in many subsequent works. This might be the reason why one frequently finds, in the complex networks literature, assertions regarding their ubiquity, their validity for almost all complex networks. Yet, the mentioned seminal works were conducted on a very limited number of networks, and, to the best of our knowledge, no large-scale study has been conducted to answer this question. In this work, we take advantage, on the one hand, of the many datasets now available online, to constitute a large collection of networks, and on the other hand, of recent analysis tools, to check the validity of this hypothesis of ubiquity. It turns out a large majority of the studied networks are indeed small-world, however this is not the case for the scale-free property, since the degree distribution of a significant proportion of our networks does not follow a power-law.}

\motscles{Réseaux complexes, propriétés topologiques, petit-monde, sans-échelle.}


\proceedings{MARAMI 2014}{-}

\maketitlepage

\section{Introduction}
\label{sec:Introduction}
Des travaux fondateurs réalisés autour des années 2000 on permis d'identifier dans les réseaux complexes des propriétés topologiques permettant de mieux comprendre le fonctionnement des systèmes modélisés. Les deux plus importantes sont certainement la propriété \textit{petit-monde}, étudiée par Watts \& Strogatz \cite{Watts1998}, et la propriété \textit{sans-échelle}, traitée par Barabási \& Albert \cite{Barabasi1999a}, qui seront toutes les deux décrites plus loin. 

À la suite de ces travaux, l'idée s'est progressivement répandue dans la littérature que ces propriétés étaient présentes dans la plupart des réseaux représentant des systèmes réels \cite{Boccaletti2006}. Pourtant, cette hypothèse d'omniprésence peut sembler fragile, pour deux raisons. Tout d'abord, les études qui se sont attachées à vérifier la prévalence de ces propriétés sont basées sur un très faible nombre de réseaux. Ainsi, les travaux fondateurs cités précédemment considèrent tous les deux seulement 3 réseaux \cite{Watts1998,Barabasi1999a}. Les publications plus récentes utilisent des collections à peine plus grandes, constituées de 5 \cite{Latora2001}, 7 \cite{Amaral2000}, 9 \cite{Boccaletti2006}, 11 \cite{Wang2003}, 15 \cite{Newman2002b}, 27 \cite{Newman2003b}, 30 \cite{Albert2002}, 33 \cite{Humphries2008} réseaux. Ceci est principalement dû au fait qu'à cette époque, ce type de données était extrêmement rare, et leur accès difficile. La seconde critique porte sur les méthodes d'analyse utilisées, en particulier pour déterminer la présence de la propriété sans-échelle, qui repose sur l'observation que le degré du réseau est distribué selon une loi de puissance. L'opération consistant à identifier cette distribution est courante en sciences expérimentales \cite{Newman2005b,Newman2011}. Des travaux méthodologiques récents ont montré qu'elle était généralement réalisée de façon inappropriée, en raison de l'utilisation d'outils statistiques inadaptés et/ou appliqués de façon incorrecte \cite{Newman2005b,Clauset2009}.

Il est aujourd'hui possible de contourner ces deux limitations. En effet, avec le développement de notre champ de recherche, des centaines, voire des milliers de réseaux différents sont maintenant accessibles via le Web, à travers les sites des chercheurs, laboratoires et journaux. Le problème principal qui se pose est de rassembler, normaliser et classifier ces données hétérogènes et éparpillées, de manière à constituer une base de données exploitable de façon automatique. Du point de vue méthodologique, là aussi, les outils permettant d'identifier la distribution du degré dans un réseau complexe existent, et la procédure à suivre a été clairement établie \cite{Clauset2009}.

Dans cet article, nous nous proposons d'éprouver l'hypothèse d'omniprésence des propriétés petit-monde et sans-échelle dans les réseaux complexes représentant des systèmes réels. Nous constituons d'abord une large collection de réseaux, puis nous l'étudions au moyen de méthodes d'analyse adéquates. Dans la section suivante, nous rappelons brièvement ce que sont ces deux propriétés qui nous intéressent. Puis, dans la section \ref{sec:Pretraitement}, nous décrivons comment nous avons constitué notre collection de réseaux, avant d'aborder en section \ref{sec:EvalProp} les méthodes d'analyse que nous avons utilisées pour l'étudier. Les résultats obtenus sont présentés dans la section \ref{sec:Resultats}. Ce travail peut-être considéré comme la première étape d'un projet à plus long terme, pour lequel nous donnons quelques pistes à explorer en guise de conclusion.

\section{Propriétés topologiques}
\label{sec:PropTopo}
De nombreuses propriétés topologiques sont aujourd'hui étudiées pour caractériser des systèmes réels, telles que la corrélation du degré, la structure de communautés ou différents types de centralité décrivant aussi bien les n\oe{}uds que les liens. Nous nous concentrons, dans cette étude, sur les deux propriétés historiques citées dans l'introduction : la propriété petit-monde et la propriété sans-échelle, que nous décrivons brièvement dans cette section.

\subsection{Petit-monde}
\label{sec:DefPetitMonde}
La propriété \textit{petit-monde} a été définie par Watts \& Strogatz en se basant sur deux mesures topologiques \cite{Watts1998}. D'une part, la distance géodésique moyenne, une caractéristique \textit{globale} du réseau qui représente le niveau de séparation typique entre deux n\oe{}uds. D'autre part, le coefficient de transitivité (ou coefficient de clustering), qui est \textit{local} et correspond au niveau d'interconnexion dans le voisinage d'un n\oe{}ud. Watts \& Strogatz ont calculé ces mesures pour différents réseaux réels, et les ont comparées à celles attendues pour les deux principaux modèles du moment : d'un côté les structures de treillis, qui constituent une forme de graphes réguliers, et qui sont qualifiés ici de complètement ordonnés ; et de l'autre les réseaux aléatoires de type Erd\H{o}s-Rényi \cite{Erdos1960}, qui sont décrits ici comme complètement désordonnés. Les premiers possèdent une transitivité et une distance moyenne élevées, alors que pour les seconds, ces deux mesures sont petites. L'analyse a montré que les réseaux réels se situent quelque part entre ces deux extrêmes. Leur distance moyenne est petite, du même ordre de grandeur que celle obtenue pour des réseaux aléatoires, et bien en-deçà de celle des treillis. Leur transitivité est élevée, du même ordre de grandeur que celle d'un treillis, et bien au-delà de celle des graphes aléatoires. C'est la présence simultanée de ces deux caractéristiques que Watts \& Strogatz ont qualifiée de propriété petit-monde. 

La propriété petit-monde a un effet sur les processus prenant place dans le système modélisé. Watts \& Strogatz ont illustré empiriquement qu'elle facilite les mécanismes de propagation dans le réseau. Par la suite, Latora \& Marchiori ont montré de façon théorique que les réseaux petit-monde sont effectivement très efficaces dans le transport d'information, à la fois aux niveaux local et global \cite{Latora2001}. Ceci a des conséquences directe sur le fonctionnement du système, qui dépendent de sa nature. Ainsi, la vitesse de propagation d'une épidémie dans une population augmente ; la capacité calculatoire d'un réseau d'automates cellulaires s'accroît ; un ensemble d'oscillateurs se synchronise plus rapidement, ou oscille de façon plus cohérente ; l'apparition d'une propriété émergente telle que la collaboration dans une population de joueurs devient plus probable ; etc. \cite{Watts1998,Lago-Fernandez2000,Li2003}. L'étude de cette propriété constitue donc un point important de l'analyse d'un système réel.






\subsection{Sans-échelle}
\label{sec:DefSansEchelle}
\vspace{-0.20cm}
Dans leur étude bien connue \cite{Barabasi1999a}, Barabási \& Albert ont observé que le degré $k$ mesuré sur plusieurs réseaux réels est distribué selon une loi de puissance $P(k) \sim k^{-\alpha}$. Autrement dit, le nombre de n\oe{}uds de degré $k$ est proportionnel à une puissance de $k$. Ils ont qualifié les réseaux concernés de \textit{sans-échelle}, en raison de la propriété qu'a cette fonction de conserver la même forme quand son paramètre est multiplié par une constante. Cependant, il faut noter que l'invariance d'échelle considérée dans cette étude ne concerne que le degré, et pas forcément les autres mesures que l'on peut calculer pour décrire un réseau (contrairement à ce que l'expression \textit{réseau sans-échelle} pourrait laisser croire). 

Comme pour la propriété petit-monde, les modèles classiques du moment n'expliquent pas cette distribution : dans les graphes réguliers, le degré est par définition le même pour tous les n\oe{}uds (loi dégénérée), alors que pour les graphes aléatoires, il est distribué selon une loi de Poisson \cite{Erdos1960}. Dans les deux cas, il existe donc une valeur caractéristique pour le degré, ce qui n'est pas le cas avec la loi de puissance. Or, la propriété sans-échelle a été observée sur des réseaux modélisant des systèmes réels très différents : grille électrique, Web, réseau d'interactions sociales, Internet, etc. \cite{Barabasi1999a,Newman2003b}. Pour l'expliquer, Barabási \& Albert utilisent deux mécanismes absents des deux modèles classiques : la croissance du système et l'attachement préférentiel. Le premier implique que le réseau soit construit itérativement, par addition de nouveaux n\oe{}uds. Le second stipule que ces nouveaux n\oe{}uds ont tendance à se connecter à des n\oe{}uds déjà existants de degré élevé. Il faut cependant noter que d'autres modèles très différents ont été proposés pour produire des réseaux sans-échelle \cite{Boccaletti2006}, ce qui laisse supposer qu'il existe plusieurs mécanismes distincts amenant à l'apparition de cette propriété dans des systèmes réels.

La propriété sans-échelle a un effet direct sur la structure du réseau : une grande majorité de n\oe{}uds est très faiblement connectée, alors qu'un nombre réduit possède un degré très élevé, et sont considérés comme des pivots (hubs). Par comparaison, dans les modèles dont la distribution du degré est plus homogène, la probabilité d'apparition de ces valeurs extrêmes est particulièrement faible. Grâce à l'existence de ces pivots, les réseaux sans-échelle opposent une meilleure résistance aux attaques aléatoires, consistant à tenter de partager le réseau en supprimant des n\oe{}uds au hasard. En revanche, ils sont plus vulnérables aux attaques concertées, qui se concentrent sur ces pivots \cite{Albert2002}. En termes de propagation, les réseaux sans-échelles modélisant des systèmes réels ont un seuil épidémiologique très bas \cite{Boccaletti2006}, ce qui signifie qu'une épidémie peut y survivre très longtemps. Comme pour la propriété petit-monde, l'effet de la propriété sans-échelle sur le fonctionnement du réseau justifie la grande importance qui lui est donnée dans la littérature.



\section{Données et prétraitement}
\label{sec:Pretraitement}
Comme nous l'avons mentionné dans l'introduction, le Web permet aujourd'hui d'accéder à une très grande quantité de réseaux complexes représentant des systèmes réels. Cependant, différents problèmes compliquent fortement la constitution d'une collection de réseaux de grande taille, qui est un prérequis à l'étude à grande échelle que nous voulons réaliser ici.

Tout d'abord, la très grande majorité des ces données librement accessibles est extrêmement dispersée, ce qui rend leur exploitation difficile. On les trouve principalement sur les pages personnelles des chercheurs qui les ont collectées ou analysées, ou bien sur les sites de certains éditeurs offrant la possibilité d'héberger des fichiers annexes que les auteurs peuvent associer à leurs articles. Certains groupes ont entrepris de constituer leurs propres collections de réseaux, ce qui permet de limiter ce problème. On peut ainsi trouver des jeux de test sur les sites de logiciels dédiés à l'analyse de réseaux, tels que Pajek\footnote{
\texttt{http://pajek.imfm.si/}
}
 ou Gephi\footnote{
\texttt{https://gephi.org/}
}. Plus récemment, des groupes de recherche ont collecté et organisé des données pour les mettre à disposition de la communautés, donnant naissance à des dépôts librement accessibles, tels que Stanford Large Network Dataset Collection\footnote{
\texttt{http://snap.stanford.edu/data/}
} ou Koblenz Network Collection\footnote{
\texttt{http://konect.uni-koblenz.de/}
}. Cependant, ces collections restent de taille relativement limitée, de l'ordre de la centaine de réseaux. 

Outre la difficulté pratique de la tâche consistant à rechercher des réseaux éparpillés sur l'ensemble du Web, des problèmes d'ordre légal empêchent aussi la constitution de grandes collections. En effet, les données ne sont pas toutes protégées par les mêmes licences d'utilisation, et l'incompatibilité de certaines d'entre elles rend impossible l'utilisation d'une licence couvrant l'ensemble de la collection. De plus, certains auteurs ou institutions refusent que leurs données soient hébergées ailleurs que sur leur propre site.

Enfin, le dernier problème est technique : il n'existe pas de format standard pour la représentation des réseaux. Ou plutôt, il en existe un grand nombre, chacun avec ces spécificités. Le format le plus simple, souvent appelé \textit{edgelist}, consiste à représenter la liste des liens du réseau, sous la forme de paires de n\oe{}uds. Cependant, ce format ne permet pas de représenter une information riche, comme par exemple des attributs nodaux, ou la dimension temporelle. Pour résoudre ce type de limitation, chaque outil propose donc son propre format, comme par exemple Gexf pour Gephi, voire plusieurs formats différents, comme c'est le cas pour Pajek. En raison de cette hétérogénéité, la constitution d'une collection requiert, en plus du travail de collecte, d'effectuer une normalisation du format de stockage des données, afin de pouvoir leur appliquer ensuite un traitement automatisé.

En raison des trois difficultés invoquées précédemment, nous avons décidé dans un premier temps de constituer une collection limitée, qui sera ensuite enrichie itérativement, et dont l'organisation sera éventuellement modifiée en fonction des retours issus de son utilisation. Lors de cette première passe, nous avons recueilli un total de $611$ réseaux, en tirant principalement parti des collections mentionnées ci-dessus, et en complétant avec des données éparses. En ce qui concerne le stockage de ces données, nous avons opté pour deux formats. Pour des réseaux que nous qualifions de \textit{basiques}, dans le sens où ils sont seulement décrits par leur structure même, avec éventuellement des liens orientés, nous avons utilisé le format \texttt{.net} de Pajek, qui est compact mais suffisant. Pour des réseaux plus riches, possédant par exemple des attributs nodaux, des liens pondérés, ou bien des types de n\oe{}uds ou de liens distincts, nous avons utilisé Graphml, qui est un dialecte XML beaucoup plus verbeux que le format de Pajek, mais dont le pouvoir expressif est en contrepartie bien plus grand. Ces deux formats offrent l'avantage d'être compatibles avec un grand nombre d'outils différents.

L'ensemble des fichiers obtenus constitue la version complète de notre collection. Cependant, pour le travail d'analyse qui nous intéresse ici, nous n'avons besoin que de réseaux basiques. Nous avons donc produit une version partielle de notre collection, contenant une version de chaque réseau débarrassée de toute information superflue. Pour la plupart des réseaux concernés, il s'agissait simplement de supprimer les orientations ou les poids des liens. Parmi les réseaux nettoyés, une centaine sont bipartis : nous avons alors calculé les deux projections possibles et gardé la moins dense (pour minimiser le fait que ces projections aboutissent à des réseaux plus denses). Il existe également une cinquantaine de réseaux multiplexes : nous avons conservé seulement un type de liens, en privilégiant celui qui permettait d'obtenir le réseau le moins dense, tout en étant connecté. Il faut souligner que ces transformations sont critiquables, car on peut argumenter qu'elles sont susceptibles d'affecter les propriétés topologiques qui nous intéressent. Nous avons néanmoins décidé d'inclure ces réseaux dans nos données, car c'est l'approche admise dans la littérature, et ce dès les travaux fondateurs précédemment mentionnés \cite{Watts1998,Barabasi1999a}. Nous projetons d'étudier plus en détail l'effet de ces transformations quand le nombre de réseaux bipartis/multiplexes dans notre collection sera suffisant. La dernière modification a consisté à supprimer les n\oe{}uds isolés (n\oe{}uds de degré nul), les liens multiples (plusieurs liens entre la même paire de n\oe{}uds) et les boucles (lien connectant un n\oe{}ud à lui-même). Enfin, nous avons également filtré les réseaux les plus grands, qu'il nous était impossible de traiter pour des raisons computationnelles. Au final, notre collection nettoyée comporte $598$ réseaux basiques, tous au format Pajek.

Pour les raisons invoquées précédemment, il nous est impossible de placer directement ces données en ligne. Cependant, la liste des réseaux concernés, avec une brève description et un lien vers les données originales, est disponible sur notre site\footnote{
\texttt{http://bit.gsu.edu.tr/compnet}
}.

\section{Évaluation des propriétés}
\label{sec:EvalProp}
Depuis les deux travaux fondateurs de la fin des années 90 mentionnés précédemment \cite{Watts1998,Barabasi1999a}, les méthodes proposées pour caractériser les propriétés petit-monde et sans-échelle ont significativement évolué, pour aller vers plus d'objectivité et de fiabilité. Nous décrivons dans cette section les méthodes que nous avons adoptées pour ce travail, et qui sont, à notre connaissance, les plus adéquates.

\subsection{Mesure petit-monde}
\label{sec:EvalPetitMonde}
Comme expliqué en section \ref{sec:DefPetitMonde}, un réseau est dit petit-monde si sa distance moyenne et sa transitivité sont respectivement du même ordre de grandeur que celles mesurées pour un réseau aléatoire et un treillis de dimension comparable. Pour mesurer cette propriété, nous avons décidé de nous baser sur la méthode proposée par Telesford \textit{et al.} \cite{Telesford2011}, qui consiste à comparer la distance moyenne $L$ et la transitivité $T$ du réseau considéré $G$, avec les valeurs $L_A$ et $T_T$ correspondant respectivement à la distance moyenne et à la transitivité mesurées pour un réseau aléatoire $G_A$ et un treillis $G_T$ considérés comme comparables à $G$. La mesure $\omega$ quantifiant la propriété petit monde est alors \cite{Telesford2011} :

\begin{equation}
\label{eq:petitmonde}
\omega(G) = \frac{L_A}{L} - \frac{T}{T_T}
\end{equation}

Pour un réseau se rapprochant d'une structure de treillis, on a $L>L_A$ et $T \simeq T_T$, donc on obtient une valeur négative. Au contraire, pour un réseau dont la distance moyenne et la transitivité ressemblent à celles d'un graphe aléatoire, on a $L \simeq L_A$ et $T < T_T$, donc la valeur de $\omega$ est positive. Enfin, si le réseau est petit-monde, alors on a simultanément $L \simeq L_A$ et $T \simeq T_T$, et $\omega$ est proche de zéro. Il faut remarquer que le fait de combiner les deux rapports pour obtenir une mesure unique n'est pas forcément pertinent, car d'autres combinaisons de $L$ et $T$ peuvent aboutir à des valeurs similaires, entraînant une interprétation erronée. Par exemple, prenons un réseau régulier correspondant à une simple grille carrée en 2D. Cette structure possède une distance moyenne très élevée, et une transitivité nulle, ce qui entraîne que les deux rapports sont quasi-nuls, et donc $\omega \simeq 0$. Or, ce réseau n'est pas du tout petit-monde, bien au contraire. Pour cette raison, il nous paraît judicieux de considérer, en plus de la mesure $\omega$, les deux termes qui la composent.

Il faut ajouter que les deux réseaux de référence $G_A$ et $G_T$ sont en fait des approximations obtenues par recâblage de $G$, via deux méthodes légèrement différentes \cite{Telesford2011}. Dans les deux cas, il s'agit de répéter un grand nombre de fois un traitement consistant à sélectionner aléatoirement deux liens et de permuter leurs n\oe{}uds de destination. La qualité de l'approximation dépend du nombre d'itérations. Notons $V$ l'ensemble des n\oe{}uds de $G$ et $E$ celui de ses liens. Considérons deux liens $(a,b)$ et $(c,d) \in E$ (où $a$, $b$, $c$ et $d \in V$) : ceux-ci deviennent $(a,d)$ et $(c,b)$ après permutation. Un filtrage adéquat des liens tirés au sort permet d'éviter de faire apparaître des liens multiples, ou de séparer le réseau en plusieurs composants distincts. Dans le cas du treillis, la permutation n'est autorisée que si elle permet de rapprocher le lien de la diagonale dans la matrice d'adjacence représentant le réseau. Autrement dit, en supposant que les n\oe{}uds sont numérotés de $1$ à $n=|V|$, on cherche à respecter la contrainte $|a-d|+|c-b| \leq |a-b|+|c-d|$. 



\subsection{Distribution du degré}
\label{sec:EvalSansEchelle}
L'identification de la loi théorique gouvernant la distribution de données empiriques est une tâche plus délicate qu'il n'y paraît. Plusieurs travaux récents on souligné le fait que les méthodes les plus répandues dans la littérature ne sont pas toujours appropriées, et ont abouti à des conclusions erronées \cite{Clauset2009}. Ceci est particulièrement vrai de la loi de puissance, qui peut être confondue avec d'autres distributions à queue lourde par certains outils traditionnellement utilisés. De plus, certaines données peuvent suivre ce type de distribution de façon partielle, i.e. seulement sur une partie de l'intervalle considéré. Dans le contexte des réseaux complexes, les caractéristiques associées à la propriété sans-échelle sont issues de la probabilité relativement élevée de voir apparaître des pivots dans le réseau. Il est donc important de s'assurer que le degré suit bien une loi de puissance, et non pas une distribution dont la queue est lourde mais décroit plus rapidement (distributions exponentielle, exponentielle étirée, log-normale, loi de puissance tronquée par une coupure exponentielle), voire qui ne possède pas du tout de queue lourde (loi de Poisson). Cette remarque vaut également si on veut valider l'hypothèse que le système étudié a été généré selon un modèle tel que celui de Barabási \& Albert, qui aboutit à la propriété sans-échelle.

Pour éviter ce type de confusion, nous utilisons la procédure proposée par Clauset \textit{et al}, qui comprend trois étapes. La première consiste à estimer les meilleures valeurs pour les paramètres $\alpha$ et $x_{min}$ de la loi, en utilisant la méthode du maximum de vraisemblance. Le premier paramètre est l'exposant de la loi, qui a déjà été mentionné en section \ref{sec:DefSansEchelle}. Le second est sa borne inférieure, qui correspond à la valeur à partir de laquelle on considère que les données analysées suivent une loi de puissance. Dans la deuxième étape, on effectue un test d'ajustement, pour vérifier s'il est plausible que les données considérées suivent la loi de puissance munie des paramètres estimés. Pour cela, on génère d'abord de nombreux jeux de données en utilisant $\alpha$ et $x_{min}$, puis on estime leurs propres paramètres. La statistique de Kolmogorov-Smirnov est ensuite utilisée pour mesurer les distances entre les distributions empiriques de ces jeux de données et les lois théoriques correspondantes (chacune étant caractérisée par les paramètres estimés). La valeur-$p$ du test d'ajustement correspond à la proportion de jeux pour lesquels on obtient une valeur de cette statistique supérieure à celle calculée pour les données originales. La troisième étape vise à comparer la loi de puissance à des distributions alternatives. Cette dernière étape est souvent ignorée dans les travaux cherchant à identifier la propriété sans-échelle, qui ne cherchent pas à vérifier si d'autres distribution sont aussi plausibles que la loi de puissance. Pourtant, même si la valeur-$p$ obtenue à la deuxième étape est suffisamment petite pour que la loi de puissance constitue une hypothèse plausible, il est possible que les données soient aussi bien, voire mieux, ajustées par d'autres distribution. On effectue la comparaison grâce au test du rapport de vraisemblance, qui consiste à comparer les vraisemblances des données étudiées, relativement aux deux distributions théoriques considérées. Le signe du logarithme de leur rapport indique laquelle des deux distributions doit être privilégiée (une valeur proche de zéro correspondant à une absence de préférence). De plus, une valeur-$p$ est calculée pour mesurer la signification de ce rapport. Le résultat final peut donc soit privilégier la loi de puissance estimée, soit identifier une ou plusieurs distributions alternatives comme plus vraisemblables, soit ne pas départager la loi de puissance et certaines alternatives.





\vspace{-0.2cm}
\section{Résultats}
\label{sec:Resultats}
\begin{figure}[thb]
{	\vspace{-0.2cm}
	\begin{center} \small 
		\includegraphics[width=0.49\textwidth]{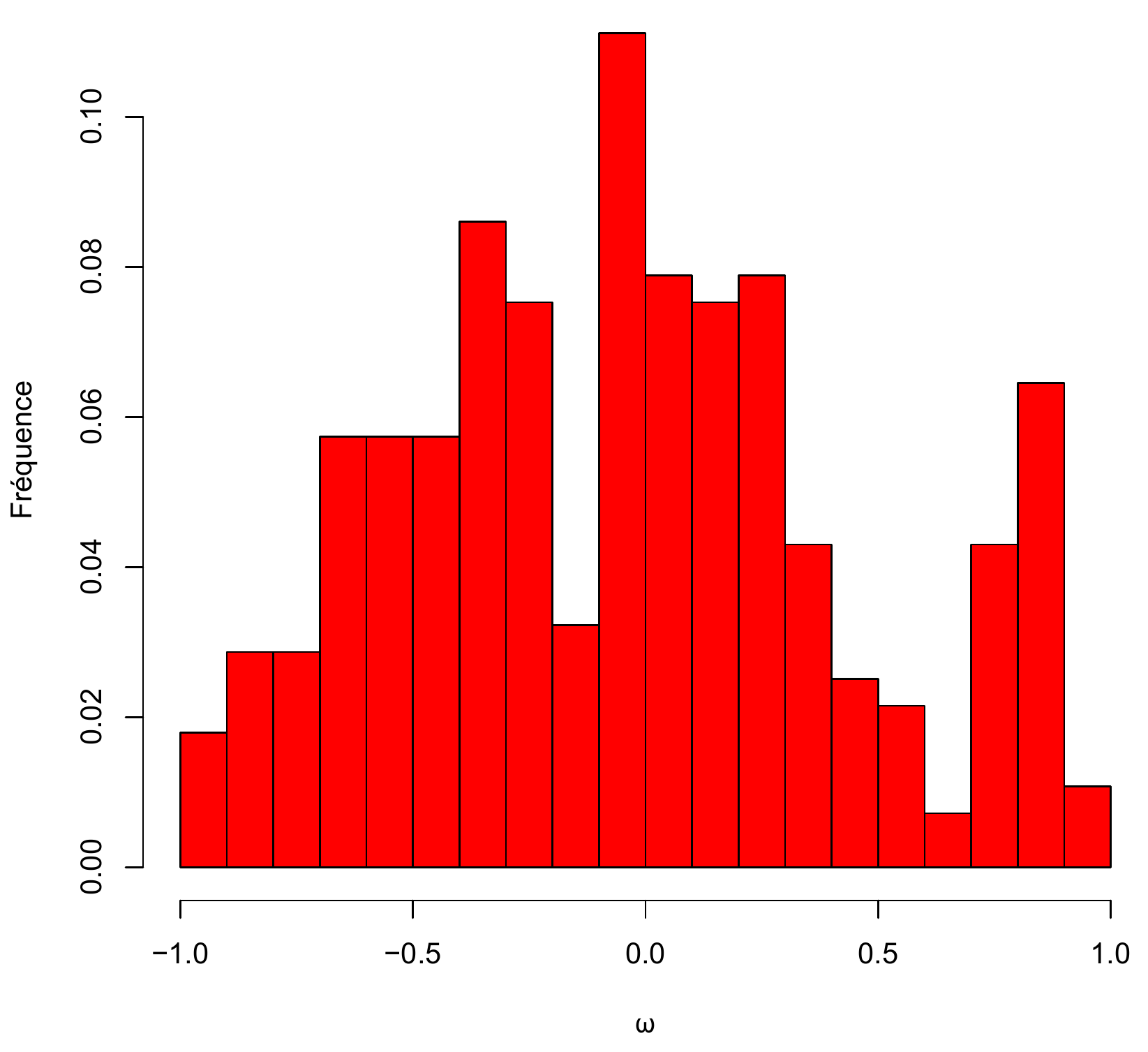}
		\includegraphics[width=0.49\textwidth]{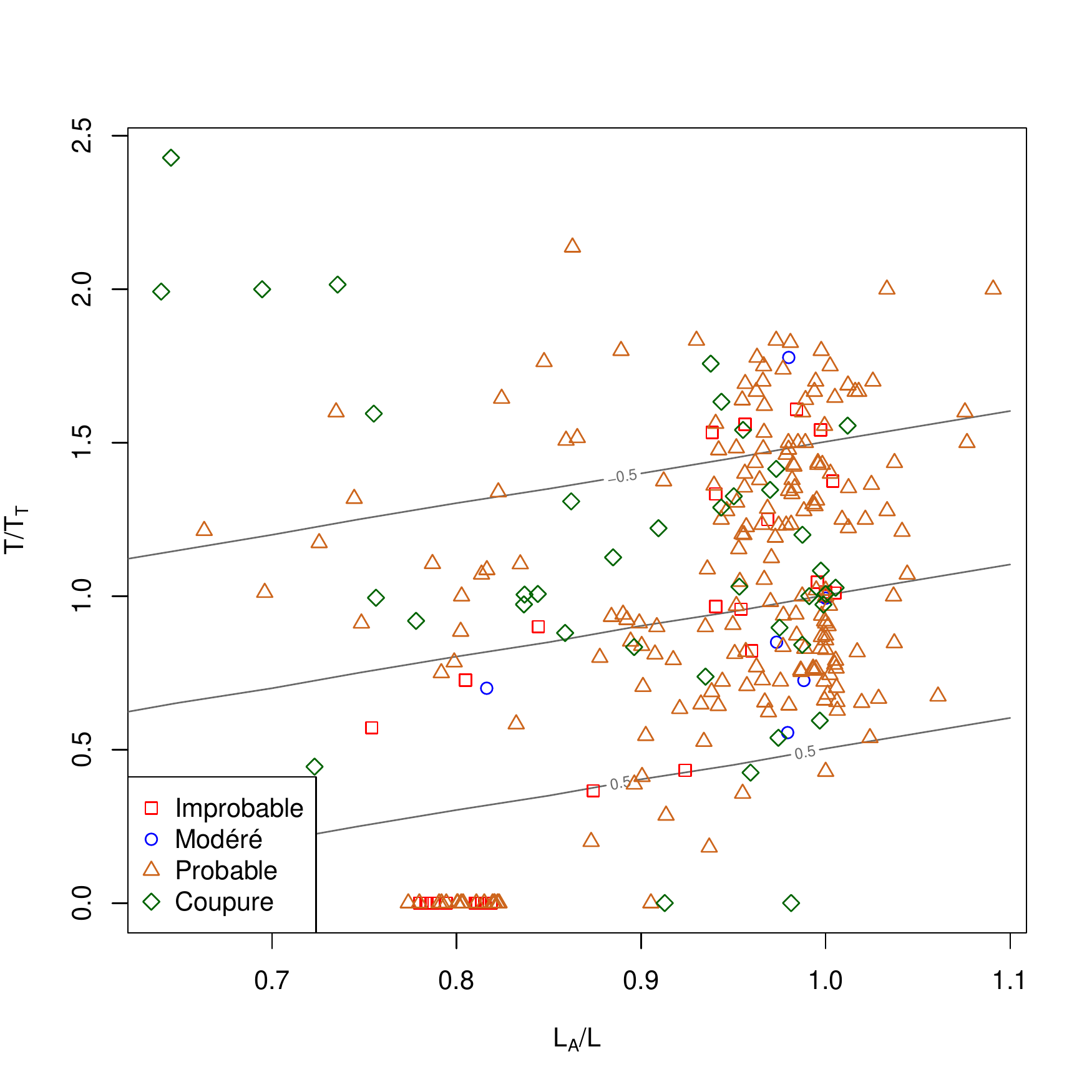}
	\end{center}\vspace{-6pt}
	\vspace{-0.5cm}
	\caption{Résultats pour la propriété petit-monde : fréquences de la mesure $\omega$ (gauche) et valeurs de ses deux composantes $L_A/L$ et $T/T_T$ (droite).}
	\vspace{-0.25cm}
	\label{fig:PetitMonde}
}
\end{figure}

Comme expliqué en section \ref{sec:EvalPetitMonde}, nous avons d'abord calculé les distances moyennes $L$ et transitivités $T$ des réseaux réels constituant notre collection, ainsi que les valeurs $L_A$ et $T_T$ de leurs versions recâblées. Avec la bibliothèque igraph\footnote{
\texttt{http://igraph.org/}
}
, on obtient pour les mesures des complexités en $O(n \langle k \rangle^2)$ pour la transitivité et $O(mn)$ pour la distance moyenne, où $n=|V|$, $m=|E|$ et $\langle k \rangle$ dénote le degré moyen du réseau. 
Le recâblage d'une paire de lien est effectué en $O(m+n)$, cependant pour obtenir de bonnes approximations, il doit être répété de très nombreuses fois, ce qui, en pratique, aboutit à des temps de calcul pouvant être relativement longs. Pour ces raisons computationnelles, il ne nous a été possible de calculer $\omega$ que pour $285$ réseaux de notre collection. L'histogramme présenté dans la Figure \ref{fig:PetitMonde} résume les résultats obtenus. Telesford \textit{et al.} \cite{Telesford2011} considèrent que les réseaux pour lequels $\omega$ se situe dans $[-0,5;0,5]$ possèdent la propriété petit-monde. Pour nos données, $65\%$ des réseaux traités sont contenus dans cet intervalle, et un test de normalité montre que la mesure $\omega$ suit approximativement une loi normale centrée sur $0$. Cela signifie notamment qu'il y a à peu près autant de réseaux penchant vers les graphes aléatoires que vers les treillis. 

Le nuage de points présenté dans la Figure \ref{fig:PetitMonde} détaille les valeurs des deux termes formant $\omega$ (cf. eq. \ref{eq:petitmonde}). Les courbes de niveaux délimitent les valeurs $\omega=-0,5$ ; $0$ et $0,5$. Les couleurs et formes représentent les résultats concernant la propriété sans-échelle, que nous discutons ci-après. On peut observer sur ce graphique qu'aucun réseau ne correspond au cas particulier mettant en défaut l'interprétation de $\omega$ (les deux termes nuls), que nous avions identifié en section \ref{sec:EvalPetitMonde}. Un groupe très dense apparait vers la coordonnée $(0,8 ;0)$, qui correspond à une collection de réseaux métaboliques présentant une transitivité nulle en raison de leur structure quasi-arborescente. Pour $72\%$ des réseaux, on a $0,9 \leq L_A/L \leq 1,1$, ce qui signifie que les réseaux tendent à avoir une distance moyenne très proche de celle de leur version randomisée. De façon générale, les réseaux ont tendance à présenter une distance moyenne plus petite ($80\%$ d'entre eux). La dispersion est plus importante pour la transitivité, puisqu'on a $0,5 \leq T/T_T \leq 1,5$ pour $67\%$ des réseaux. Ceci indique donc une tendance plus importante à s'éloigner de la transitivité obtenue sur leurs versions randomisées. On peut aussi remarquer qu'à la différence de la distance, les écarts observés sont répartis de façon relativement similaire des deux côtés de la valeur $1$.

Les tests concernant la distribution du degré sont moins complexes, calculatoirement parlant, et ont pu être menés sur $510$ réseaux. Nous avons utilisé la nomenclature de Clauset \textit{et al}. \cite{Clauset2009}, qui distingue les situations suivantes : \textit{Improbable} (la distribution ne suit probablement pas une loi de puissance), \textit{Modéré} (la loi de puissance ajuste bien les données, mais c'est aussi le cas d'autres distributions), \textit{Probable} (la loi de puissance est la seule à bien ajuster les données) et \textit{Coupure} (loi de puissance tronquée par une coupure exponentielle est la meilleure alternative). La loi de puissance est la seule alternative probable pour $55\%$ des réseaux testés (niveau Probable), alors que pour $1\%$ elle est une alternative parmi d'autres (niveau Modéré). Pour $28\%$ des réseaux testés, le degré suit probablement une loi de puissance tronquée (niveau Coupure), et les réseaux qui ne suivent probablement pas une loi de puissance représentent $16\%$ (niveau Improbable). Au total, on a donc $44\%$ de réseaux dont le degré peut être caractérisé par une échelle type, et qui ne possèdent donc pas la propriété sans-échelle. Ce résultat est en contradiction avec l'hypothèse d'omniprésence de cette propriété parmi les réseaux complexes.

Le nuage de points inclus dans la Figure \ref{fig:PetitMonde} indique les résultats relatifs à la distribution du degré à l'aide d'un code de couleur. On peut remarquer que les réseaux de niveau Improbable semblent groupés dans la zone du graphique comprise entre $\omega=-0,5$ et $0.5$, mais ceci pourrait être dû au faible nombre de réseaux concernés. Pour les niveaux Probable et Coupure (rappelons que la propriété sans-échelle n'est pas valide pour le second), la distribution est sensiblement similaire, avec la même concentration autour du point $(1;1)$. Cela semble indiquer que, sur ces données, les propriétés petit-monde et sans-échelles ne sont pas liées.

\begin{figure}[thb]
{	\vspace{-0.2cm}
	\begin{center} \small 
		\includegraphics[width=0.49\textwidth]{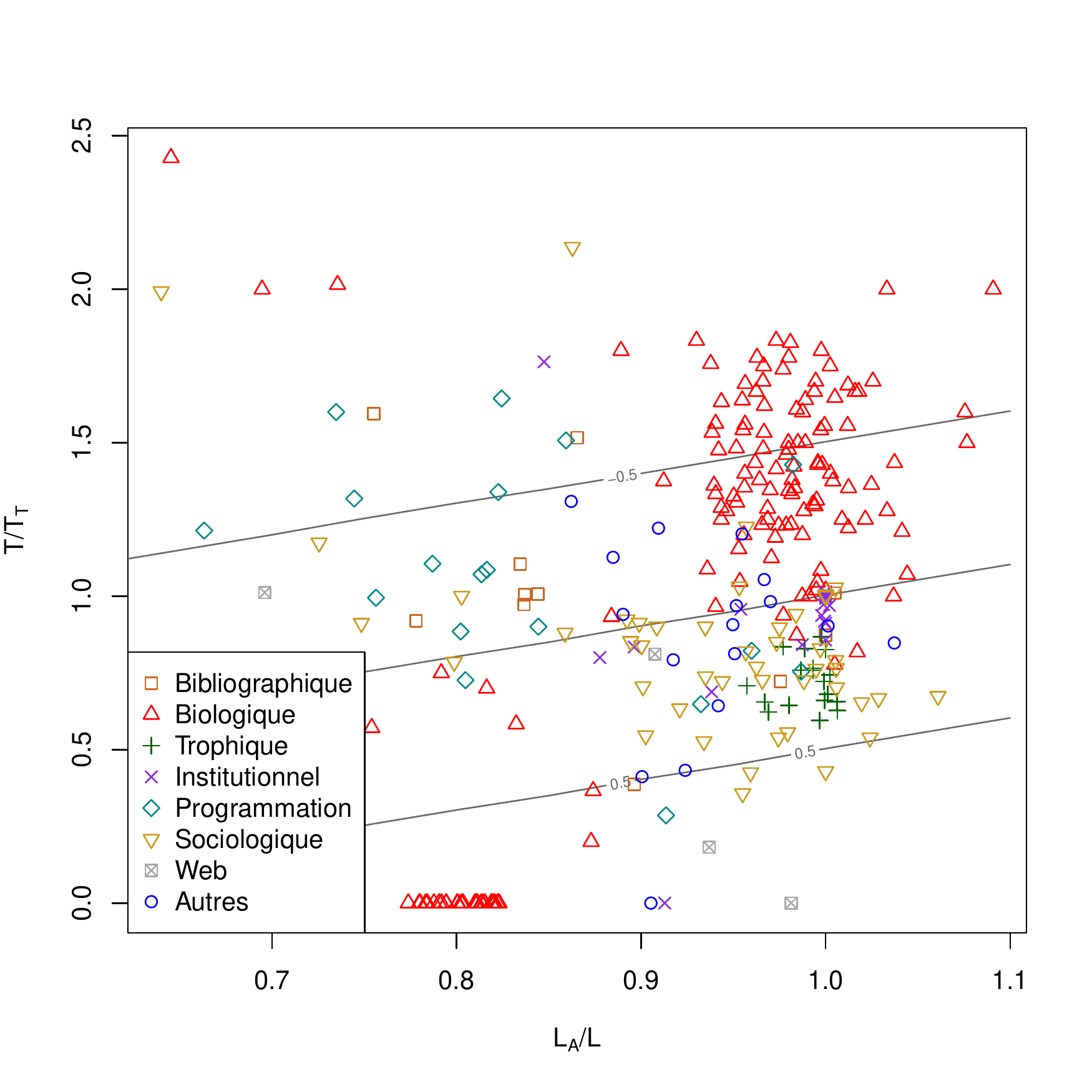}
	\end{center}\vspace{-6pt}
	\vspace{-0.5cm}
	\caption{Distribution des classes de systèmes.}
	\label{fig:Classes}
	\vspace{-0.3cm}
}
\end{figure}

L'analyse de la relation entre la nature du système modélisé et la présence (ou l'absence) des deux propriétés qui nous intéressent constitue un autre point intéressant. Dans les articles portant sur l'analyse d'un nombre suffisamment important de réseaux, il est fréquent de voir les auteurs tenter de faire des regroupements, afin de généraliser leurs résultats. Ces groupes sont généralement basés sur les systèmes modélisés. Par exemple, les réseaux ayant trait au fonctionnement des être vivants sont regroupés sous le nom de réseaux biologiques, les services de réseautage en ligne sous le nom de réseaux sociaux, etc. \cite{Albert2002,Newman2003b}. Nous avons fait de même avec nos données, et la Figure \ref{fig:Classes} donne une représentation graphique des principales classes définies. On peut voir que certaines occupent des zones bien spécifiques de l'espace : réseaux biologiques, réseaux d'interactions sociales, réseaux trophiques. Cependant, ce n'est pas vrai pour d'autres : réseaux bibliographiques, institutionnels et de programmes. La question demande à être creusée, et nous pensons nous orienter vers des classes de relations plutôt que de systèmes : interaction, similarité, référence, appartenance, etc.

\section{Conclusion}
\label{sec:Conclusion}
Dans le travail présenté ici, nous nous somme concentrés sur la constitution d'une collection de réseaux, et sur l'application de méthodes récentes afin d'éprouver l'hypothèse d'omniprésence des propriétés petit-monde et sans-échelle, qui apparait très fréquemment dans la littérature. Nos résultats, pour l'instant partiels, montrent que si la propriété petit-monde est effectivement présente dans une large majorité de réseaux, avec toutefois des fluctuations non négligeables en termes de distance moyenne et de transitivité, en revanche près de la moitié des réseaux testés ne possède pas la propriété sans-échelle.

Le but de cet article était de présenter la première étape d'un projet à plus long terme, visant à analyser de façon systématique une grande quantité de réseaux complexes réels. L'idée est de généraliser les travaux fondateurs de Watts \& Strogatz \cite{Watts1998} et de Barabási \& Albert \cite{Barabasi1999a}, qui se concentraient sur un petit nombre de réseaux et de mesures topologiques. La quantité de données librement accessibles est beaucoup plus importante aujourd'hui, et de nombreuses mesures ont été proposées depuis, aussi nous nous proposons d'exploiter ces ressources pour mettre à jour de nouvelles régularités du même type que les propriétés petit-monde et sans-échelle. La suite directe de se travail consistera à explorer la relation entre ces deux propriétés : sous quelles conditions les réseaux sans-échelle sont-ils aussi petit-monde ? À plus longue échéance, nous souhaitons étendre notre étude en termes de données et d'outils. D'abord en complétant notre collection de réseaux, puis en calculant les autres principales mesures topologiques existantes : centralité, modularité, etc. Il sera vraisemblablement nécessaire d'appliquer des méthodes d'analyse et de fouille de données pour détecter des régularités, soit valables pour toute la collection, soit caractéristiques de certaines classes de systèmes.

\vspace{-0.3cm}
\bibliography{marami14}

\end{document}